%% file: shifted-GMRES-recyc.tex
\documentclass[final]{siamltex}
\usepackage{euscript,graphicx,epsfig,amsfonts,mathtools,epstopdf}
\usepackage[ruled,lined,vlined,linesnumbered,algosection]{algorithm2e}

\input macros
\def\RGMRES{Recycled GMRES}

\title{Krylov Subspace Recycling for Sequences of Shifted Linear Systems\thanks{%
This version dated \today. This research was supported in part by the U.S. Department of Energy under grant DE-FG02-05ER25672 and the U.S.\ National Science Foundation under grant DMS-1115520.}}

\author{Kirk M. Soodhalter\footnotemark[2],\ \ 
	Daniel B. Szyld\footnotemark[3],\ \ 
        and Fei Xue\footnotemark[4]}

\begin{document}
\maketitle
\renewcommand{\thefootnote}{\fnsymbol{footnote}} 
\footnotetext[2]{Industrial Mathematics Institute, Johannes Kepler University, 
Altenbergerstra{\ss}e 69, A-4040 Linz, Austria.
({\tt soodhalter@indmath.uni-linz.ac.at})}
\footnotetext[3]{Department of Mathematics, Temple University,
1805 N Broad Street, Philadelphia, PA 19122-6094.
({\tt szyld@temple.edu}).}
\footnotetext[4]{Department of Mathematics, University of Louisiana at Lafayette, Lafayette, LA 70504.
({\tt fxue@louisiana.edu}) }

\begin{abstract}
We study the use of Krylov subspace recycling for the solution of a sequence of 
slowly-changing families of linear systems, where each family consists of shifted 
linear systems that differ in the coefficient matrix only by multiples of the identity. 
Our aim is to explore the simultaneous solution of each family of shifted systems within 
the framework of subspace recycling, using one augmented subspace to extract candidate 
solutions for all the shifted systems.  The ideal method would 
use the same augmented subspace for all systems and have fixed storage requirements, independent
of the number of shifted systems per family.  We show that
a method satisfying both requirements cannot exist in this framework.  

As an alternative, we introduce two schemes.  One constructs a separate deflation space 
for each shifted system but solves
each family of shifted systems simultaneously.  The other builds only one recycled subspace and
constructs approximate corrections to the solutions of the shifted systems at each cycle of the iterative
linear solver while only minimizing the 
base system residual.  At convergence of the base system solution, we apply the method recursively to the remaining 
unconverged systems. 
We present numerical examples involving systems arising 
in lattice quantum chromodynamics.
\end{abstract}

\begin{keywords}
Krylov subspace methods, subspace recycling, shifted linear systems, QCD
\end{keywords}


\pagestyle{myheadings}
\thispagestyle{plain}
\markboth{K. M. SOODHALTER, D. B. SZYLD, AND F. XUE}{SUBSPACE RECYCLING FOR SHIFTED LINEAR SYSTEMS}

\input intro.tex

\input prelim.tex

\input collinear-nonexist.tex

\input shifted-GCRODR.tex

\input approx-collinear-correction.tex

\input numerical_results.tex

\input conclusions.tex

\input acknowledgement.tex

\bibliographystyle{siam}
\bibliography{../../../master}
\end{document}

%% file: macros.tex
\def\cal{\mathcal}
\def\norm#1{\left\|#1\right\|}

\def\C{\mathbb{C}}

\def\R{\mathbb{R}}

\def\Cnn{\C^{n\times n}}

\def\Rn{\R^n}

\def\BA{{\bf A}}  
  
  \def\CC{{\cal C}}

  \def\CF{{\cal F}}

  \def\CK{{\cal K}}

  \def\CO{{\cal O}}
  \def\CP{{\cal P}}
  
  \def\CR{{\cal R}}
  \def\CS{{\cal S}}
  
  \def\CU{{\cal U}}

\def\BAs{\BA{\kern-1.5pt}}

\def\CPs{\CP{\kern-0.8pt}}

\mathcode`\@="8000
{\catcode`\@=\active \gdef@{\mkern1mu}}

\def\mydate{\number\day\ {\ifcase\month \or January\or February\or
              March\or April\or May\or June\or July\or August\or
              September\or October\or November\or December\fi}
\number\year}
\def\vek#1{\mathbf{#1}}

\newcommand{\argmin}[1]{\underset{#1}{\text{{\rm argmin}}}}
\newcommand{\floor}[1]{\lfloor #1 \rfloor}

\newcommand{\prn}[1]{\left(#1\right)}

\newcommand{\brac}[1]{\left[#1\right]}
\newcommand{\curl}[1]{\left\{#1\right\}}
\newcommand{\ab}[1]{\left|#1\right|}

\def\cvd{~\vbox{\hrule\hbox{%
  \vrule height1.3ex\hskip0.8ex\vrule}\hrule } }

%% file: intro.tex
\section{Introduction}\label{intro}
We consider the solution of a sequence of families of non-Hermitian linear systems.  Let $\CF$ 
denote a family of coefficient matrices differing by multiples of the identity.  In other words, 
\begin{equation}\label{eqn.shifted-family}
\CF = \curl{\vek A+\sigma^{(\ell)}\vek I}_{\ell=1}^{L}\subset \Cnn,
\end{equation}
where $L$ is the number of matrices in the family, and we are solving the family of linear systems
\begin{equation}\label{eqn.one-shifted-fam}
	\prn{\vek A+\sigma^{(\ell)}\vek I}\vek x^{(\ell)}=\vek b\mbox{\ \ \ for\ \ \ }\ell=1,\ldots , L.
\end{equation}

We call the numbers $\curl{\sigma^{(\ell)}}_{\ell=1}^{L}\subset\C$ \textit{shifts}, $\vek A$ 
the \textit{base matrix}, and $\vek A+\sigma\vek I$ a \textit{shifted matrix}.  
Systems of the form (\ref{eqn.one-shifted-fam}) are called \textit{shifted linear systems}.  
There are many applications which warrant the solution of a family of shifted linear 
systems with coefficient matrices belonging to $\CF$,  such as  those arising in lattice quantum 
chromodynamics (QCD) 
(see, e.g., \cite{Frommer1995}) as well as other applications such as Tikhonov-Philips regularization, 
global methods of nonlinear analysis, and Newton trust region 
methods \cite{CDT1985}.  Krylov subspace methods have been 
proposed to simultaneously solve this family of systems 
\cite{Frommer2003}, \cite{Frommer1998}, \cite{Simoncini2003a}.  

Our goal is to explore simultaneously solving a family of shifted systems (or a sequence of families)
over an augmented Krylov subspace, i.e., we explore how one would incorporate 
existing shifted system techniques into the subspace recycling framework 
\cite{Parks.deSturler.GCRODR.2005}.  Does a method exist which (a) simultaneously 
solves all systems in the family, using one subspace to extract candidate solutions, 
and (b) satisfies a fixed storage requirement, independent of 
the number of shifts?

In this paper, we treat this question in the context of the recycled GMRES framework 
\cite{Parks.deSturler.GCRODR.2005} combined with GMRES for shifted systems \cite{Frommer1998}.
We demonstrate that the two mentioned requirements (a) and (b) cannot be 
achieved simultaneously.  
We present two methods: one which sacrifices fixed storage and the
other which sacrifices the simultaneous solution of all shifted systems in each family, 
instead solving one system 
and simultaneously improving the approximations of the others at a very modest cost.

For simplicity (and to avoid excessive indices), our discussion will mostly center around
solving the model problem,
\begin{eqnarray}
\vek A\vek x & = &\vek b \label{eqn.Axb}\\
 (\vek A+\sigma\vek I)\vek x^{(\sigma)} & = &\vek b, \label{eqn.Axb-shifted}
\end{eqnarray}
using \RGMRES\ in the presence of a $k$-dimensional initial recycled subspace~$\CU$.
However, what we derive applies to a more general situation.
Let $\CF_{i}$ denote  the $i$th family of linear systems, defined by 
\begin{equation}\nonumber
\CF_{i} = \curl{\vek A_{i}+\sigma_{i}^{(\ell)}\vek I}_{\ell=1}^{L_{i}}\subset\Cnn, 
\end{equation}
where $L_{i}$ denotes the number of linear systems to be solved at step $i$.  In other words, at step $i$, for shifts $\curl{\sigma_{i}^{(\ell)}}_{\ell=1}^{L_{i}}$ we are 
solving systems of the form
\[
	\prn{\vek A_{i}+\sigma_{i}^{(\ell)}\vek I}\vek x_{i}^{(\ell)} = \vek b_{i}\mbox{\ \ \ for\ \ \ }\ell=1\ldots L_{i}.
\]
This is the general problem our desired method is meant to address.

It should be noted that, for solving the model problem in the absence of the initial subspace $\CU$, 
techniques already have been developed
to solve a family of systems simultaneously, building deflation subspaces 
from harmonic Ritz vectors; see, e.g., \cite{Darnell2008}.  We are exploring here the situation
in which we have an initial deflation subspace $\CU$, a case for which the techniques
presented in \cite{Darnell2008} do not apply.

In the next section, we review some existing methods for solving (\ref{eqn.Axb}) and (\ref{eqn.Axb-shifted}), and 
we describe the framework of subspace recycling used in, e.g., \cite{Parks.deSturler.GCRODR.2005}.  
In Section \ref{section.nonexistence-collin}, we show that it is generally \textit{not} possible to construct 
solutions for all shifted systems over the same augmented subspace in a way that is compatible with restarting
while allowing for the simultaneous solution of all systems. 
In Section \ref{sect.shifted-gcrodr}, we present a method which sacrifices the fixed storage requirement.  
This method is a direct extension of the one presented in \cite{Frommer1998}.    
We present a scheme in Section \ref{approx-collinear} which sacrifices the requirement that 
each family be solved simultaneously.
This method produces improved approximations for the 
shifted system while solving the base system, using only one recycled subspace.  This method is 
also derived from \cite{Frommer1998}, but the approximations are not computed
according to the same residual collinearity constraints as described in \cite{Frommer1998}.  
In Section \ref{numresults}, we 
present numerical results for a family of simple bidiagonal matrices and for some sequences 
of QCD matrices obtained from \cite{DH.2011} and \cite{Wuppertal.Matrices.2012}.

%% file: prelim.tex
\section{Preliminaries}\label{sect.prelim}
In many Krylov subspace iterative methods, recall that we generate an orthonormal basis for 
the Krylov subspace
\[ 
\CK_{j}(\vek A,\vek u) = \text{span}\curl{\vek u, \vek A\vek u,\ldots, \vek A^{j-1}\vek u}
\] 
with the Arnoldi process, where $\vek u$ is some starting vector.  Let $\vek V_{j}\in\C^{n\times j}$ be the 
matrix with orthonormal columns generated by the Arnoldi process spanning $\CK_{j}(\vek A,\vek u)$.  
Then we have the Arnoldi relation
\begin{equation}\label{eqn.arnoldi-relation}
	\vek A\vek V_{j} = \vek V_{j+1}\overline{\vek H}_{j}
\end{equation}
with $\overline{\vek H}_{j}\in\C^{(j+1)\times j}$; see, e.g., 
\cite[Section 6.3]{Saad.Iter.Meth.Sparse.2003} and 
\cite{szyld.simoncini.survey.2007}.  Let $\vek x_{0}$ be an initial approximation and 
$\vek r_{0}=\vek b-\vek A\vek x_{0}$ be the initial residual. At iteration $j$, 
we compute $\vek x_{j} = \vek x_{0} + \vek t_{j}$, 
where $\vek t_{j}\in\CK_{j}(\vek A, \vek r_{0})$.  
In GMRES \cite{Saad.GMRES.1986}, we choose
\[
	\vek t_{j} = \argmin{\vek t\in \CK_{j}(\vek A,\vek r_{0})}\norm{\vek b-\vek A(\vek x_{0} + \vek t)},
\]
and this is equivalent to solving the smaller minimization problem
\begin{equation}\label{eqn.GMRES-least-squares}
	\vek y_{j} = \argmin{\vek y\in\C^{j}}\norm{ \overline{\vek H}_{j}\vek y - \norm{\vek r_{0}}\vek e_{1}^{(j+1)}},
\end{equation}
where we use the notation $\vek e_{\ell}^{(k)}$ to denote the $\ell$th Cartesian basis vector in $\R^{k}$,
and setting $\vek x_{j} = \vek x_{0} + \vek V_{j}\vek y_{j}$.  In restarted GMRES (GMRES($m$)), we halt this process 
at step $m$, discard the matrix $\vek V_{m}$, and restart with the new initial residual
\linebreak$\vek r_{0} \leftarrow \vek b - \vek A \vek x_{m}$.  
This process is repeated until we achieve convergence.  Adaptions of
 restarted GMRES to solve (\ref{eqn.one-shifted-fam}) have been previously proposed; 
 see, e.g., \cite{Frommer1998}.  

It should be noted that methods based on the nonsymmetric Lanczos process have also been adapted 
for solving (\ref{eqn.one-shifted-fam}).  
Extensions of methods, such as BiCGStab~\cite{Frommer2003} and QMR, have been 
developed \cite{Frommer1995}.  
A recently proposed method called IDR \cite{Sonneveld2008}, which has been shown to be a 
generalization of BiCGStab \cite{Simoncini2010a}, has also been extended to 
solve~(\ref{eqn.Axb})~and (\ref{eqn.Axb-shifted}) \cite{Kirchner2011}. We will not deal with 
nonsymmetric Lanczos-based methods in this paper, but these alternatives are worth mentioning.

Many methods for solving (\ref{eqn.one-shifted-fam}) use the fact that for any 
shift $\sigma$, the Krylov subspace generated by $\vek A$ and $\vek b$ is invariant under the shift, i.e., 
\begin{equation}\nonumber
	\CK_{j}(\vek A,\vek b) = \CK_{j}(\vek A+\sigma\vek I,\widetilde{\vek b}),
\end{equation} 
as long as the starting vectors are collinear, i.e., $\widetilde{\vek b}=\beta\vek b$, with a shifted Arnoldi 
relation similar to~(\ref{eqn.arnoldi-relation})
\begin{equation}\label{eqn.shifted-arnoldi-relation}
	(\vek A+\sigma\vek I)\vek V_{j} = \vek V_{j+1}\overline{\vek H}_{j}^{(\sigma)}.
\end{equation}
Note that the shift-invariance no longer holds if general preconditioning is used. 
However, polynomial preconditioning \cite{Jegerlehner.shifted-systems.1996} would be appropriate in this setting.  
There has been recent work on choosing optimal polynomial preconditioners in the setting of solving multiple shifted 
systems \cite{ASV.2012}, \cite{M.2013}.  In this project, though, we focus on the 
unpreconditioned case, as in \cite{Frommer1998}, \cite{Morgan.GMRESDR.2002}.

The shift-invariance property indicates that large savings in storage and 
time can be achieved by generating only one sequence of Krylov subspaces and solving all shifted 
systems in one Krylov subspace simultaneously.  Suppose that the initial residuals of (\ref{eqn.Axb}) and 
(\ref{eqn.Axb-shifted}) are collinear.  As we iterate, we simply apply the Petrov-Galerkin condition
with the same subspaces for 
all residuals.  However, once restarting is introduced, the situation becomes more complicated.  
The projected residuals may no longer be collinear and the Krylov subspaces at restart will not be equivalent.
In \cite{Frommer2003}, a general theorem is presented which
describes conditions under which the residuals will be naturally collinear in this manner.
In \cite{Frommer1998} it is observed that
the GMRES residual projection does not have this property.   

Frommer and Gl\"{a}ssner \cite{Frommer1998} proposed a restarted GMRES method to solve 
(\ref{eqn.Axb})--(\ref{eqn.Axb-shifted}).  Suppose that 
the residuals for the shifted and base systems are collinear, i.e., 
$\vek r_{0}^{(\sigma)}=\beta_{0}\vek r_{0}$. 
Within a cycle,
for the base system approximation, the residual is minimized using GMRES.  The shifted system
approximation is found by requiring the residual to be collinear to that of the base system, i.e., 
\begin{equation}\label{eqn.mth-resid-colin}
	\vek r^{(\sigma)} _{m} = \beta_{m}\vek r_{m}.
\end{equation}
After computing the GMRES solution for the base system, we can denote the GMRES
least-squares residual $\vek z_{m+1} = \norm{\vek r_{0}}\vek e^{(m+1)}_{1} 
- \overline{\vek H}_{m}\vek y_{m}$.
It is shown in \cite{Frommer1998} that 
 for (\ref{eqn.mth-resid-colin}) to hold, we must have
\[ 
	\overline{\vek H}^{(\sigma)}_{m}\vek y^{(\sigma)} _{m} + \vek z_{m+1}\beta_{m} = \beta_{0}\norm{\vek r_{0}}\vek e_{1}^{(m+1)},
\] 
and $\vek x^{(\sigma)} _{m} = \vek x^{(\sigma)} _{0} + \vek V_{m}\vek y^{(\sigma)} _{m}$,
where 
\begin{equation}\nonumber
	\overline{\vek H}_{m}^{(\sigma)} = \overline{\vek H}_{m} + \begin{bmatrix}  \sigma\vek I_{m\times m} \\ \vek 0_{1\times m}\end{bmatrix}.
\end{equation}
Thus, we can compute both $\vek y^{(\sigma)} _{m}$ and $\beta_{m}$ by solving the augmented linear system,
\begin{equation}\label{eqn.mth-resid-equation}
	\brac{\begin{matrix}\overline{\vek H}^{(\sigma)}_{m} &  \vek z_{m+1} \end{matrix}}\brac{\begin{matrix}\vek y^{(\sigma)} _{m}\\ \beta_{m}\end{matrix}} = 
\beta_{0}\norm{\vek r_{0}}\vek e_{1}^{(m+1)}.
\end{equation}

The collinear residual exists if and only if the residual polynomial $r_{m}(t)$, 
associated with $\vek r_{m}$ satisfies 
$r_{m}(-\sigma)\neq 0$; otherwise, the augmented system is singular  \cite[Lemmas 2.1 and 2.4]{Frommer1998}.  
 For a positive-real matrix $\vek A$ (field of values being contained in the right half-plane), 
restarted GMRES for shifted linear systems computes solutions at every iteration for all shifts 
$\sigma^{(i)} > 0$ and, in addition, we have 
$\norm{\vek r_{m}} \leq \norm{\vek r_{m}^{(\sigma_{i})}}$ for such shifts \cite{Frommer1998}.  
The shifts applied in the setting of QCD yield a family of coefficient 
matrices which are, in theory, real-positive \cite{Frommer1998}.  

We briefly review the \RGMRES\  method described in \cite{Parks.deSturler.GCRODR.2005}.
This algorithm represents the confluence of two approaches: those descending from 
the implicitly restarted Arnoldi method \cite{LS.1996}, 
such as Morgan's GMRES-DR \cite{Morgan.GMRESDR.2002}, 
and those descending from de Sturler's GCRO method \cite{deSturler.GCRO.1996}. 
GMRES-DR is a restarted GMRES algorithm, where at the end of each cycle, harmonic Ritz vectors
are computed, and a subset of them are used to augment the Krylov subspace generated at the next cycle.
The GCRO method allows the user to select the optimal correction over arbitrary subspaces.  
This concept is extended by de Sturler in \cite{deSturler.GCROT.1999}, where 
a framework is provided for selecting the optimal 
subspace to retain from one cycle to the next so as to minimize the error produced by 
discarding useful information accumulated in the subspace for candidate solutions before restart.  
This algorithm is called GCROT, where OT stands for optimal truncation.  A simplified version of the GCROT approach, 
based on restarted GMRES (called LGMRES) is presented in \cite{Manteuffel.LGMRES.2005}. 
Parks et~al.\ in
\cite{Parks.deSturler.GCRODR.2005} combine the ideas of  \cite{Morgan.GMRESDR.2002} 
and \cite{deSturler.GCROT.1999} 
and extend them to a sequence of slowly-changing linear systems. They call their method GCRO-DR (\RGMRES). 

Suppose we are solving (\ref{eqn.Axb}), and we have a $k$-dimensional subspace $\CU$
whose image under the action of $\vek A$ is $\CC = \vek A\,\CU$.  Let $\vek{P}$ be the 
orthogonal projector onto $\CC^{\perp}$.  Furthermore, let $\vek x_{0}$ be such that
$\vek r_{0}\in\CC^{\perp}$ (this is always cheaply available).  We generate the Krylov subspace with respect to the 
projected operator $\vek{P}\vek A$, $\CK_{m}(\vek{P}\vek A,\vek r_{0})$.  At iteration $m$, 
the \RGMRES\ method generates the approximation 
\begin{equation}\nonumber
	\vek x_{m}=\vek x_{0}+\vek s_{m} + \vek t_{m}
\end{equation}
where $\vek s_{m}\in\CU$ and $\vek t_{m}\in\CK_{m}(\vek{P}\vek A,\vek r_{0})$.  The corrections
$\vek s_{m}$ and $\vek t_{m}$ are chosen according to the minimum residual, Petrov-Galerkin
condition over the augmented Krylov subspace, i.e.,
\begin{equation}\label{eqn.rgmres-petrov-galerkin}
	\vek r_{m}\perp \vek A\prn{\CU + \CK_{m}\prn{\vek{P}\vek A,\vek r_{0}}}.  
\end{equation}
At the end of the cycle,
an updated $\CU$ is constructed, the Krylov subspace basis is discarded, and we restart.
At convergence, $\CU$ is saved, to be used when solving the next linear system.

In terms of implementation, \RGMRES\ can be described as a modification of the GMRES
algorithm.  Let $\vek U\in\C^{n\times k}$ have columns spanning $\CU$, scaled such that
$ 
	\vek C=\vek A\vek U
$ 
has orthonormal columns.  Then we can explicitly construct
$\vek{P} = \vek I - \vek C\vek C^{\ast}$.  At each iteration, applying $\vek{P}$ is equivalent
to performing $k$ steps of the Modified Gram-Schmidt process to orthogonalize
the new Arnoldi vector against the columns of $\vek C$. The orthogonalization coefficients
generated at step $m$ are stored in the $m$th column of $\vek B_{m} = \vek C^{\ast}\vek A\vek V_{m}$,
and $\vek B_{m+1}$ is simply $\vek B_{m}$ with one new column appended.
Let $\overline{\vek H}_{m}$ and $\vek V_{m}$ be defined as before, but for the projected
Krylov subspace $\CK_{m}\prn{\vek{P}\vek A,\vek r_{0}}$.  Enforcing (\ref{eqn.rgmres-petrov-galerkin})
is equivalent to solving the GMRES minimization problem (\ref{eqn.GMRES-least-squares}) for
 $\CK_{m}\prn{\vek{P}\vek A,\vek r_{0}}$ and setting 
 \begin{equation}\nonumber
	\vek s_{m}=-\vek U\vek B_{m}\vek y_{m}\mbox{\ \ and\ \ }\vek t_{m}=\vek V_{m}\vek y_{m},
\end{equation}
so that 
\begin{align} \nonumber
	\vek x_{m}& = \vek x_{0} - \vek U\vek B_{m}\vek y_{m} + \vek V_{m}\vek y_{m} 
= \vek x_{0} + \begin{bmatrix} \vek U & \vek V_{m}  \end{bmatrix}
			       \begin{bmatrix}  -\vek B_{m}\vek y_{m}\\ \vek y_{m}\end{bmatrix} .
\end{align}
This is a consequence of the fact that the \RGMRES\ least squares problem, as stated in \cite[Equation 2.13]{Parks.deSturler.GCRODR.2005} can be satisfied exactly in the first $k$ rows.

Convergence results for augmented Krylov subspace methods were shown in,
e.g., \cite{Eiermann.Analysis-accel.2000, Saad.Deflated-Aug-Krylov.1997}, but not much work 
has been done in the context of \RGMRES .  Some not-yet-published work has been presented by 
de Sturler that specifically addresses the convergence behavior of optimal methods in which we recycle 
using the above framework \cite{desturler.recycling-convergence-abstract.2012}.  This work asserts that 
the improvement of convergence bounds from recycling a particular subspace can be quantified according 
to the quality of the recycled subspace as an invariant subspace of $\vek A$.  A particular finding, backed up by 
empirical observation, is that an approximate invariant subspace of modest quality (as judged by the 
largest principal angle between $\vek U$ and $\vek C$) will still yield improvements in bounds on the 
residual norm.

It should be noted that for a single system, that deflation and seeding of the Krylov subspace in the context
of shifted systems (and specifically QCD) have been previously considered; see, e.g., 
\cite{AMW.2008},\cite{SO.2010}.
Furthermore, if we have no initial recycled space and compute harmonic Ritz vectors at each restart, 
\RGMRES\ is algebraically equivalent to Morgan's GMRES-DR \cite{Morgan.GMRESDR.2002}.  
Iterating orthogonally to an approximate invariant subspace to accelerate convergence of GMRES 
can be justified by the theoretical work in \cite{Simoncini2005}.  It was shown that the widely observed 
two-stage convergence behavior of GMRES, which has been termed \textit{superlinear convergence}, 
is governed by how well the Krylov subspace approximates a certain eigenspace.  Specifically, when 
the Krylov subspace contains a good approximation to the eigenspace (call this eigenspace $\CS$) 
associated to eigenvalues hindering convergence, we will switch from the slow phase to the fast phase, 
and convergence will mimic that of GMRES on the projected operator $\vek{P}_{\CS}^{\perp}\vek A$ where 
$\vek{P}_{\CS}^{\perp}$ is the orthogonal projector onto the orthogonal complement of $\CS$.  This 
analysis complements previous discussions of this phenomenon, 
see e.g., \cite{Ipsen.GMRES-minimal-poly.1996}, \cite{Vorst1993a}.

%% file: collinear-nonexist.tex
\section{Nonexistence of the Ideal Method}\label{section.nonexistence-collin}

Subspace recycling has shown great potential to improve the convergence of restarted methods, in many 
cases, without dramatically increasing memory costs.  Therefore, if we can 
incorporate GMRES for shifted linear systems into the recycling framework described 
in \cite{Parks.deSturler.GCRODR.2005}, we will have a storage-efficient method which 
will solve all shifted systems simultaneously. 
%
In this context, it is most natural to consider extending GMRES for 
shifted systems \cite{Frommer1998} into the recycling framework.
We denote such a 
method  \textit{\RGMRES\ for shifted systems}.  
We explore how such an algorithm would look and show that 
we generally cannot satisfy the fixed memory 
requirement while achieving simultaneous solution of all systems using a single augmented subspace.
%

Consider the simplified model problem, with linear systems (\ref{eqn.Axb})--(\ref{eqn.Axb-shifted}), 
subspaces $\CU$ and $\CC$, and their respective matrix counterparts $\vek U$ and $\vek C$.
The ideal method will solve (\ref{eqn.Axb}) using \RGMRES\ while generating 
approximations for (\ref{eqn.Axb-shifted}) of the form
\begin{equation}\label{xst.eq}
	\vek x_{m}^{(\sigma)} = \vek x_{0}^{(\sigma)} + \vek s_{m}^{(\sigma)} + \vek t_{m}^{(\sigma)}
\end{equation}
with $\vek s_{m}^{(\sigma)}\in\CU$ and 
$\vek t_{m}^{(\sigma)}\in \CK_{m}(\vek P\vek A,\vek r_{0})$, such that we have residual collinearity.  
Such a method could be used for any number of shifts without increasing
storage requirements.

We begin with a useful result about Krylov subspaces for projected operators.
\begin{proposition}\label{prop.proj-shift-kryl-eqiv}
	Let $\vek C$ be a matrix with orthonormal columns spanning $\CC$.  
	Then $\vek v\in\CC^{\perp}$, i.e., $\vek C\vek C^{\ast}\vek v=\vek 
0$, if and only if
\begin{equation}\label{eqn.prop.proj-shift-kryl-eqiv}
	 \CK_{m}(\vek P\vek A,\vek v) = \CK_{m}(\vek P(\vek A+\sigma\vek I),\vek v)\mbox{\ for all\ } m
\end{equation}
\end{proposition}
\textit{Proof.} 
First, suppose $\vek v\perp \CC$. Since $\vek C\vek C^{\ast}\vek v = \vek 0$, we have
\[
	(\vek I-\vek C\vek C^{\ast})(\vek A+\sigma\vek I)\vek v = (\vek I-\vek C\vek C^{\ast})\vek A\vek v + \sigma(\vek I-\vek C\vek C^{\ast})\vek v = (\vek I-\vek C\vek 
C^{\ast})\vek A\vek v + \sigma\vek v.
\]
Therefore, when restricted to vectors orthogonal to  $\CR(\vek C)$, we have that 
\[
    (\vek I-\vek C\vek C^{\ast})(\vek A+\sigma\vek I) = (\vek I-\vek C\vek C^{\ast})\vek A + \sigma\vek I.
\]
Furthermore, since any $\vek u\in\CR(\vek P(\vek A+\sigma\vek I))$ is orthogonal to $\CR(\vek C)$, we have
\[
		\brac{(\vek I-\vek C\vek C^{\ast})(\vek A+\sigma\vek I)}^{j}\vek v = \brac{(\vek I-\vek C\vek C^{\ast})\vek A + \sigma\vek I}^{j}\vek v
\]
when applied to any $\vek v\perp \CR(\vek C)$. Thus,
\[
\CK_{m}(\vek P(\vek A + \sigma\vek I),\vek v) = \CK_{m}(\vek P\vek A + \sigma\vek I,\vek v) = \CK_{m}(\vek P\vek A,\vek v),
\]
where the last equality follows from the shift invariance property of Krylov subspaces.

\def\Ac{(\vek I-\vek C\vek C^{\ast})\vek A}
Conversely, suppose (\ref{eqn.prop.proj-shift-kryl-eqiv}) holds, and let $m=2$.  
Due to the equivalence of the two subspaces, for any vector  $\vek u\in\CK_{2}(\vek P\vek A,\vek v)\setminus\CK_{1}(\vek P\vek A,\vek v)$ we have
\begin{equation}\nonumber
\vek u = \alpha_{1}\vek v + \alpha_{2}\Ac\vek v = \beta_{1}\vek v + \beta_{2}\Ac\vek v + \beta_{2}\sigma(\vek I-\vek C\vek C^{\ast})\vek v
\end{equation}
where $\alpha_{2}$ and $\beta_{2}$ are nonzero.  This implies
\begin{equation}\nonumber
(\alpha_{2}-\beta_{2})\Ac\vek v - \beta_{2}\sigma(\vek I-\vek C\vek C^{\ast})\vek v = (\beta_{1}-\alpha_{1})\vek v,
\end{equation}
and thus, $\vek v\perp\CR(\vek C)$. ~~\cvd

Thus, for $\vek r_{0}\in\CC^{\perp}$, the projected Krylov subspace is invariant 
under a constant shift of the matrix $\vek A$, and the shifted Arnoldi relation (\ref{eqn.shifted-arnoldi-relation})
holds as well.  

In \cite{Parks.deSturler.GCRODR.2005}, it is shown that the augmented Krylov subspace satisfies
an Arnoldi-like relation, namely
\begin{equation}\label{eqn.GCRODR.arnoldi}
	\vek A\widehat{\vek V}_{m} = \widehat{\vek W}_{m+1}\overline{\vek G}_{m},
\end{equation}
where 
\begin{equation}\nonumber 
	\widehat{\vek V}_{m} = \brac{\begin{matrix}{\vek U} & \vek V_{m}\end{matrix}}\mbox{,\ }\widehat{\vek W}_{m+1} = \brac{\begin{matrix}\vek C & 
\vek V_{m+1}\end{matrix}}\mbox{, and\ \ }\overline{\vek G}_{m} = \brac{\begin{matrix}\vek I_{k} & \vek B_{m}\\ \vek 0 & \overline{\vek H}_{m}\end{matrix}}.
\end{equation}
Even with Proposition \ref{prop.proj-shift-kryl-eqiv}, the relation (\ref{eqn.GCRODR.arnoldi}) does not have
a shifted analog, as in (\ref{eqn.shifted-arnoldi-relation}).  
Instead, we have
\[
	(\vek A + \sigma\vek I)\widehat{\vek V}_{m}  = \widehat{\vek W}_{m+1} \brac{\begin{matrix}\vek I_{k} & \vek B_{m}\\ \vek 0 & \overline{\vek H}_{m} \end{matrix}} + \sigma \widehat{\vek V}_{m}.
\]
If we have 
\begin{equation}\label{eqn.range-containment}
	\CR(\widehat{\vek V}_{m}) \subset \CR(\widehat{\vek W}_{m+1}),
\end{equation} 
the relation could be easily modified so that a relation similar to (\ref{eqn.shifted-arnoldi-relation}) 
holds, allowing the collinearity condition to be enforced.  
However, this inclusion, in general, does not hold; the columns of $\vek U$ 
might span an approximate invariant subspace of $\vek A$, not a true invariant subspace.  
Similar observations are made in the context of Hermitian systems in 
\cite{Kilmer.deSturler.tomography.2006}.     

There is at least one scenario in which (\ref{eqn.range-containment}) does hold.  
Consider the situation in which we begin with no starting recycled space and compute 
harmonic Ritz vectors at the end of each cycle to pass to the next cycle.  
We run an $m$-step cycle of shifted GMRES, and at the end of that cycle, let the columns
 of $\vek U$ be $k$ 
harmonic Ritz vectors, we compute $\vek C$ as before, and restart.  
Morgan \cite{Morgan2000} showed  that for a harmonic Ritz pair $(\vek g,\theta)$, the 
eigenvector residual $\vek A\vek g - \theta\vek g$ is a multiple of the GMRES 
residual $\vek r_{m}$.  At the end of a cycle, if we compute $k$ harmonic Ritz vectors 
and store them as the columns of $\widetilde{\vek U}$, then we know that
\begin{equation}\label{eqn.AU-UD.span}
\CR(\vek A\widetilde{\vek U} - \widetilde{\vek U}\vek D) = \text{span}(\vek r_{m}),
\end{equation}
 where $\vek D = \text{diag}(\theta_{1},\ldots,\theta_{k})$, the diagonal matrix containing 
 the harmonic Ritz values associated to the columns of $
 \widetilde{\vek U}$.  If we compute the QR-factorization 
 of $\vek A\widetilde{\vek U}=\vek C\vek R$ and let $\vek U = \widetilde{\vek U}\vek R^{-1}$, then for $\vek T 
 = \vek R\vek D\vek R^{-1}$ we have 
\begin{equation}\nonumber
	\CR(\vek C - \vek U\vek T) = \text{span}(\vek r_{m}).  
\end{equation}
 At the beginning of the next cycle, we take 
$ 
\vek v_{1} = \vek r_{m}/\norm{\vek r_{m}}
$ 
 as the first Krylov vector; and in this case, the containment (\ref{eqn.range-containment}) holds.  
This is the same fact exploited in \cite{Darnell2008}, where the authors observe that the augmented Krylov subspace is itself actually a larger
Krylov subspace with a different starting vector.  Thus, the shifted GMRES method can be applied directly to the Krylov subspace augmented
with the harmonic Ritz vectors, as long as there was no deflation space at the beginning of the process.     

What about in the general setting?  Let $\vek E$ be a matrix whose columns form a basis for the orthogonal complement of $\CC~\oplus~\CK_{m+1}(\vek P\vek A,\vek r_{0})$ in $\Rn$.  We note that $\vek E$ needs not be computed; we use it here as a theoretical tool.  We can write 
\begin{equation}\label{eqn.U-decompose}
	\vek U = \vek C\vek Y + \vek V_{m+1}\vek Z + \vek E\vek F,
\end{equation}
where $\vek Y\in\C^{k\times k}$, $\vek Z\in\C^{(m+1)\times k}$, and $\vek F\in \C^{(n-m-1-k)\times k}$.  This yields the following \textit{imperfect} Arnoldi-like relation for the shifted system,
\begin{equation}\label{eqn.gcro-shifted-arnoldi}
	(\vek A + \sigma\vek I) \brac{\begin{matrix}\vek U & \vek V_{m}\end{matrix}} = \brac{\begin{matrix}\vek C & \vek V_{m+1}\end{matrix}}\brac{\begin{matrix}\vek I_{k} + \sigma\vek Y & \vek B\\ \vek \sigma\vek Z & \overline{\vek H}_{m}^{(\sigma)} \end{matrix}} + \sigma \brac{\begin{matrix}\vek E\vek F & \vek 0\end{matrix}}.
\end{equation}
If we let 
\begin{equation}\nonumber 
	\widetilde{\vek G}_{m}^{(\sigma)} = \brac{\begin{matrix}\vek I_{k} + \sigma\vek Y & \vek B\\ \vek \sigma\vek Z & \overline{\vek H}_{m}^{(\sigma)} \end{matrix}},
\end{equation}
together with (\ref{eqn.GCRODR.arnoldi}), then the Arnoldi-like relation (\ref{eqn.gcro-shifted-arnoldi}) can be rewritten as
\[
	(\vek A + \sigma\vek I)\widehat{\vek V}_{m} = \widehat{\vek W}_{m+1}\widetilde{\vek G}_{m}^{(\sigma)} + \sigma \brac{\begin{matrix}\vek E\vek F & \vek 0\end{matrix}}.
\]
We can write the correction $\vek s_{m}$ and $\vek t_{m}$ obtained by the \RGMRES\ minimization as,
\begin{equation} \label{st.eq}
	\vek s_{m} = \vek U\vek y^{(1)}_{m}\mbox{\ \ and\ \ }\vek t_{m}=\vek V_{m}\vek y^{(2)}_{m},
\end{equation}
and stack $\vek y^{(1)}_{m}$ and $\vek y^{(2)}_{m}$ in the vector 
\begin{equation}\nonumber
	\widehat{\vek y}_{m}=\begin{bmatrix} \vek y^{(1)}_{m}\\ \vek y^{(2)}_{m} \end{bmatrix}.
\end{equation}
In \cite{Parks.deSturler.GCRODR.2005}, the \RGMRES\ minimization is written so that we are
computing $\widehat{\vek y}_{m}$, satisfying
 \begin{eqnarray}
  	\vek r_{m}  & = & \vek r_{0} - \widehat{\vek W}_{m+1}\overline{\vek G}_{m}\widehat{\vek y}_{m} \nonumber \\
  	& = & \norm{\vek r_{0}}\widehat{\vek W}_{m+1}\vek e_{k+1}^{(m+1)} - \widehat{\vek W}_{m+1}\overline{\vek G}_{m}\widehat{\vek y}_{m} = \widehat{\vek W}_{m+1}\widehat{\vek z}_{m+1},\nonumber
 \end{eqnarray} 
where we used (\ref{eqn.GCRODR.arnoldi}) and 
\begin{equation} \label{eqn.z-def}
\widehat{\vek z}_{m+1} = \norm{\vek r_{0}}\vek e^{(m+1)}_{k+1} - \overline{\vek G}_{m}\widehat{\vek y}_{m}  
\end{equation}
is the \RGMRES\ least-squares residual.
Now, for the shifted system, we would like to enforce the collinearity condition.  If a collinear residual were to exist for the shifted system, then it would satisfy
\begin{align}
 	\vek r_{m}^{(\sigma)} & =  \beta_{m}\vek r _{m} &\iff\nonumber \\
 	\vek b - (\vek A + \sigma\vek I)(\vek x_{0}^{(\sigma)} + \widehat{\vek V}_{m}\vek y^{(\sigma)} _{m}) & =  \beta_{m}\widehat{\vek W}_{m+1}\widehat{\vek z}_{m+1}  &\iff\nonumber \\
	\vek r_{0}^{(\sigma)} - (\vek A + \sigma\vek I)\widehat{\vek V}_{m}\vek y^{(\sigma)} _{m} & =  \widehat{\vek W}_{m+1}\widehat{\vek z}_{m+1}\beta_{m}  &\iff\nonumber \\
	\beta_{0}\vek r_{0} - (\widehat{\vek W}_{m+1}\widetilde{\vek G}_{m}^{(\sigma)} + \sigma \brac{\begin{matrix}\vek E\vek F & \vek 0\end{matrix}})\widehat{\vek y}^{(\sigma)} _{m} 
& =  \widehat{\vek W}_{m+1}\widehat{\vek z}_{m+1}\beta_{m}  &\iff\nonumber \\
	\beta_{0}\vek r_{0} & =  \widehat{\vek W}_{m+1}(\widehat{\vek z}_{m+1}\beta_{m} + \widetilde{\vek G}_{m}^{(\sigma)}\widehat{\vek y}^{(\sigma)} _{m})\nonumber \\
	                                &\quad\, + \sigma \brac{\begin{matrix}\vek E\vek F & \vek 0\end{matrix}}\widehat{\vek y}_{m}^{(\sigma)}. \label{eqn.gcrodr-colin-resid}
\end{align}
Observe that in the general case, $\vek r_{0}\in\CC \oplus \CK_{m}(\vek P\vek A,\vek r_{0})$ 
while the right-hand side of  (\ref{eqn.gcrodr-colin-resid}) has 
a non-zero component in $\CR(\vek E) = (\CC \oplus  \CK_{m}(\vek P\vek A,\vek r_{0}))^{\perp}$.   
Thus, we state the conditions for 
existence (and nonexistence) of the collinear residual in the following theorem.
\begin{theorem}\label{thm.collinear-resid-exist}
	Suppose we have approximations $\vek x_{0}$ and $\vek x_{0}^{(\sigma)}$ 
	to the solutions of {\rm (\ref{eqn.Axb})} and {\rm (\ref{eqn.Axb-shifted})}, respectively, such that the residuals $\vek r_{0}$ and $\vek r_{0}^{(\sigma)}$ 
	are collinear, and $\vek r_{0}\in\CC^{\perp}$.  Let $\vek r_{m}$ be the minimum residual solution 
	produced by \RGMRES\  over the augmented Krylov subspace $\CU+\CK_{m}(\vek P\vek A,\vek r_{0})$.  
	Then one of the following is true:
	\begin{itemize}
		\item $\CU + \CK_{m}(\vek P\vek A,\vek r_{0}) \subset \CC \oplus \CK_{m+1}(\vek P\vek A,\vek r_{0}) $
		\item There exists \textbf{no} approximation 
		$\vek x_{m}^{(\sigma)}\in\CU+\CK_{m}(\vek P\vek A,\vek r_{0})$ to 
		{\rm (\ref{eqn.Axb-shifted})} such that $\vek r_{m}^{(\sigma)}$ is collinear to 
		$\vek r_{m}$, i.e., $\vek r_{m}^{(\sigma)} \neq \beta_{m}\vek r_{m}$,\mbox{\ for \textbf{all}\ } $\beta_{m}\in\C$.
	\end{itemize}  
\end{theorem}

%% file: shifted-GCRODR.tex
\section{A Method with the Colinearity Approach}
\label{sect.shifted-gcrodr}
We have shown in Theorem~\ref{thm.collinear-resid-exist}
that the ideal algorithm, i.e., one where all shifted systems
are solved with the same approximation subspace and with fixed storage,
generally does not exist.  However, by removing
one of the two requirements, we can derive viable methods.  
First, in this section, we consider a method which imposes the
collinearity of the residuals, thus allowing the use of the
same subspace for all shifts, at the cost of building
different deflation subspaces for each of the shifts with their dimensions small enough to
not incur excessive memory costs.
In other words, additional storage is required for each new shift.
We do so by extending the work of Frommer and Gl\"{a}ssner~\cite{Frommer1998}
to this situation.

Let ${\vek x}_{-1}$ and ${\vek x}^{(\sigma)}_{-1}$ be initial approximations
so that the initial residuals are 
collinear, i.e., ${\vek r}_{-1} = {\beta}_{0}{\vek r}^{(\sigma)}_{-1}$.  The update,
\begin{equation}\label{eqn.init-proj-update}
	\vek x_{0} = \vek x_{-1} + \vek U\vek C^{\ast}\vek r_{-1}\mbox{\ \ and\ \ }\vek r_{0} = \vek r_{-1} - \vek C\vek C^{\ast}\vek r_{-1}
\end{equation}
cheaply yields a residual $\vek r_{0}\in\CC^{\perp}$.
In order to effect a similar update of $\vek x_{-1}^{(\sigma)}$, 
we need $\vek U^{(\sigma)}$ such that 
\begin{equation}\label{eqn.AUCsig}
	\vek C = \vek A\vek U = (\vek A+\sigma \vek I)\vek U^{(\sigma)}. 
\end{equation}
This requires an additional $k$ vectors of storage for each shift.
Given an initial subspace $\CU$, we can derive $\vek U^{(\sigma)}$ for each value of $\sigma$.
For details, 
see \cite{SSX.2013}, where in addition to a description of how to efficiently build 
the family deflation spaces, an analysis is presented on the relation between
the value of the shift and the degradation of the orthogonality of
the columns of the matrix in (\ref{eqn.AUCsig}).

As described in Section~\ref{section.nonexistence-collin}, (\ref{eqn.Axb}) can be solved using \RGMRES\ while 
for (\ref{eqn.Axb-shifted}), we can compute 
\begin{equation}\label{eqn.shifted-update} \nonumber
	\vek x_{m}^{(\sigma)} = \vek x_{0}^{(\sigma)} + \vek U^{(\sigma)}\vek y_{m}^{(1,\sigma)} + \vek V_{m}\vek y_{m}^{(2,\sigma)}
\end{equation}
such that the collinearity condition $\vek r_{m}^{(\sigma)} = \beta_{m}^{(\sigma)} \vek r_{m}$ holds;
cf.~(\ref{st.eq}). The vector $\vek y_{m}^{(\sigma)}$  together with the scalar
$\beta_m^{(\sigma)}$ 
are solved simultaneously from an augmented system, as was done 
in~\cite{Frommer1998}. Just as \RGMRES\ can be viewed as applying GMRES to a 
projected linear system, this method can be shown to reduce to applying the shifted GMRES
method to a projected, shifted linear system.

Such an augmented system is also used in our
second approach presented in the next section; cf.~(\ref{eqn.approx-collinear-system}). 
The procedure for this method with multiple deflation 
spaces is fully developed in \cite{SSX.2013},
but omitted here for sake of brevity.
We observe though, that the approximation with collinear residual is drawn
from $\vek U^{(\sigma)} + \CK_{m}(\vek P\vek A,\vek r_{0})$ rather than from 
$\vek U + \CK_{m}(\vek P\vek A,\vek r_{0})$, from which the minimal residual correction 
of the base system is extracted.

%% file: approx-collinear-correction.tex
\section{A Method with Fixed Storage}\label{approx-collinear}

Inspired by the results of Theorem~\ref{thm.collinear-resid-exist}, we consider a different
alternative than that briefly discussed in Section~\ref{sect.shifted-gcrodr}.
If we enforce the fixed-storage requirement (i.e., only one recycled subspace $\CU$ is stored
and all approximations are drawn from the same augmented Krylov subspace) then 
a prospective algorithm must overcome two obstacles. 

First, we cannot conveniently update the residual of the shifted system.  
For the shifted system, we construct
approximations of the form (\ref{xst.eq}).
As already discussed, without a $\vek U^{(\sigma)}$ defined as in (\ref{eqn.AUCsig}), 
we cannot project $\vek r_{-1}^{(\sigma)}$ 
and update $\vek x_{-1}^{(\sigma)}$, as in (\ref{eqn.init-proj-update}).  
 As a remedy,
 we can perform an update of the shifted system approximation which implicitly updates the 
 residual by the perturbation of an orthogonal projection.
 We set 
 \begin{equation}\nonumber
	\vek x_{0}^{(\sigma)} = {\vek x}_{-1}^{(\sigma)} + \vek U\vek C^{\ast}{\vek r}_{-1}^{(\sigma)}.
\end{equation}
The updated residual can be written as
\begin{eqnarray}
 {\vek r}_{0}^{(\sigma)}& = &  \vek b-(\vek A+\sigma\vek I)\vek x_{0}^{(\sigma)}\nonumber\\
			 	  & = &  \vek b-(\vek A+\sigma\vek I)({\vek x}_{-1}^{(\sigma)}+\vek U\vek C^{\ast}{\vek r}_{-1}^{(\sigma)}) \nonumber\\
			 	  & = &  {\vek r}_{-1}^{(\sigma)}-(\vek A+\sigma\vek I)\vek U\vek C^{\ast}{\vek r}_{-1}^{(\sigma)} \nonumber\\
			 	  & = &  \underbrace{{\vek r}_{-1}^{(\sigma)}- \vek C\vek C^{\ast}{\vek r}_{-1}^{(\sigma)}}_{\text{true orthogonal projection}} -\underbrace{\sigma\vek U\vek C^{\ast}{\vek r}_{-1}^{(\sigma)} . }_{perturbation}\label{eqn.proj-pert}
\end{eqnarray}

Second, the collinear residual does not exist. 
Deriving this result yields clues to another way forward.  
 Neglecting a term from (\ref{eqn.gcrodr-colin-resid}) allows us to solve a nearby approximate collinearity 
 condition (which we will explain shortly, after Algorithm \ref{alg.GCRODR-approx-collinear}) and update the 
 approximation for the shifted system.  This update is of the 
 form (\ref{xst.eq}) with $\vek s_{m}^{(\sigma)}\in\CU$.
 These corrections tend to improve the residual but do not lead to convergence for the 
 shifted system, which will start with an expected improved approximation. 
 We present analysis showing how much improvement is possible with this method.  
 After convergence of the base system, 
 the algorithm can be applied recursively on the remaining unconverged systems.  
 This recursive method of solving one seed system at a time while choosing corrections 
 for the approximations for the other systems has been previously suggested in the context 
 of linear systems with multiple
 right-hand sides; see e.g., \cite{Chan1997},\cite{SPM1.1981}.

We begin by providing
an overview of the strategy we are proposing and encode this into a schematic algorithm.  
This algorithm solves the base system with \RGMRES\ while cheaply computing better initial 
approximations for the shifted systems.
We present this outline in Algorithm \ref{alg.GCRODR-approx-collinear}.

\begin{algorithm}[hb!]
\caption{Schematic of Shifted \RGMRES\  with an Approximate Collinearity Condition}
\label{alg.GCRODR-approx-collinear}
\SetKwInOut{Input}{Input}\SetKwInOut{Output}{Output}\SetKwComment{Comment}{}{}
\Input{$\vek A\in\C^{n\times n}$; $\curl{\sigma^{(\ell)}}_{\ell=1}^{L}\subset\C$; $\vek U,\vek C\in\C^{n\times k}$ such that $\vek A\vek U=\vek C$ and $\vek C^{\ast}\vek C = \vek I_{k}$; Initial Approximations $\vek x_{0}$ and $\vek x_{0}^{(\sigma^{(\ell)})}$ such that residuals are collinear; $\varepsilon > 0$}
	$\vek x \leftarrow \vek x_{0}$, $\vek r = \vek b-\vek A\vek x$\\
	$\vek x\leftarrow \vek x + \vek U\vek C^{\ast}\vek r$, $\vek r \leftarrow \vek r - \vek C\vek C^{\ast}\vek r$\Comment*[r]{Project base residual}
	$\vek x^{(\sigma^{(\ell)})}\leftarrow \vek x_{0}^{(\sigma^{(\ell)})}$, $\vek r ^{(\sigma^{(\ell)})}= \vek b-\vek A\vek x^{(\sigma^{(\ell)})}$ for all $\ell$\\
	\For{$\ell=1\text{ to }L$}{
		$\vek x^{(\sigma^{(\ell)})}\leftarrow \vek x^{(\sigma^{(\ell)})} + \vek U\vek C^{\ast}\vek r^{(\sigma^{(\ell)})}$\Comment*[r]{Update shifted approximation, but  not an implicit residual projection}\label{alg.line.noproj}
	}
	\While{$\norm{\vek r}>\vek \varepsilon$}{
		Construct a basis of the subspace $\CK_{m}((\vek I-\vek C\vek C^{\ast})\vek A,\vek r)$\\
		Compute update $\vek t\in \cal{R}(\vek U) + \CK_{m}((\vek I-\vek C\vek C^{\ast})\vek A,\vek r)$ by minimizing residual using 
		\RGMRES \\
		$\vek x\leftarrow\vek x+\vek t$; $\vek r\leftarrow\vek b-\vek A\vek x$\\
		\For{$\ell=1\text{ to }L$}{
			Compute update $\vek t^{(\sigma^{(\ell)})}\in \cal{R}(\vek U) + \CK_{m}((\vek I-\vek C\vek C^{\ast})\vek A,\vek r)$ according to the approximate collinearity condition\label{alg.line.approx-collin}\\
			$\vek x^{(\sigma^{(\ell)})} \leftarrow \vek x^{(\sigma^{(\ell)})} + \vek t^{(\sigma^{(\ell)})}$
		}
	
		Compute updated recycled subspace information $\vek U$ and $\vek C$
	}
	Clear any variables no longer needed\\
	\If{$L>2$}{
		Make a recursive call to Algorithm \ref{alg.GCRODR-approx-collinear} with $\vek A\leftarrow \vek A + \sigma^{(1)}\vek I$, shifts $\curl{\sigma^{(\ell)}-\sigma^{(1)}}_{\ell=2}^{L}$, 
approximations $\curl{\vek x^{(\sigma^{(\ell)})}}_{\ell=2}^{L}$ and updated recycled subspace matrix $\vek U$
	}
	\Else{
		Apply \RGMRES\ to the last unconverged system
	}
\end{algorithm}

This algorithm relies on dropping the term $\vek E\vek F\prn{\tilde{\vek y}_{m}}_{1:k}$
from (\ref{eqn.gcrodr-colin-resid}),
which yields an augmented linear system that can be solved directly,
\begin{eqnarray}
 \vek z_{m+1}\tilde{\beta}_{m} + \widetilde{\vek G}_{m}^{(\sigma)}\tilde{\vek y}^{(\sigma)} _{m}& = & \beta_{0}\norm{\vek r_0}\vek e_{k+1}^{(m+1)}\label{eqn.approx-collinear-system}\mbox{,\ \ \ or}\\
\brac{\begin{matrix}  \widetilde{\vek G}_{m}^{(\sigma)} &   \vek z_{m+1} \end{matrix}}\brac{\begin{matrix}\tilde{\vek y}^{(\sigma)} _{m}\\ \tilde{\beta}_{m} \end{matrix}} & = &  \beta_{0}\norm{\vek r_0}\vek e_{k+1}^{(m+1)}.\nonumber
\end{eqnarray}
Thus, we proceed by solving this nearby problem and updating the shifted solution,
\begin{equation}\label{eqn.shifted-incorrect-update}
	\vek x_{m}^{(\sigma)} = \vek x_{0}^{(\sigma)} +\widehat{\vek V}_{m}\tilde{\vek y}^{(\sigma)} _{m}.
\end{equation}
For each restart cycle, we repeat this process for the shifted system.  We stop when the base residual 
norm is below tolerance.  
When the residual norm for the base 
system reaches the desired tolerance, the residual norm 
of the shifted system will have been reduced at little additional cost;
but, generally, the reduction is insufficient.  
Thus, we apply the GMRES with 
recycling algorithm with this approximate collinearity scheme to the remaining unsolved systems, 
taking one of the shifted systems as
 our new base system.  This method is amenable to recursion on the number of shifts.  
 When only one system remains, 
\RGMRES\ is applied.  

Observe that for any number of shifts, we can easily form $\widetilde{\vek G}_{m}^{(\sigma)}$ 
for each $\sigma$ at little additional cost.  
The matrices $\vek Y$ and $\vek Z$ in (\ref{eqn.U-decompose}) must be computed only once 
per cycle, regardless of the number of shifted 
systems we are solving.  However, additional shifts will require more recursive calls to the 
algorithm and, thus, more iterations.


Why does the approximate collinearity condition 
produce an improved approximation to the solution of the shifted system? 
How well we can expect the algorithm to perform?  The following analysis answers these
questions and also 
yields a cheap way in which we can monitor the progress of the residuals of the shifted systems.
Theorem~\ref{thm.gcrodr-shift-convergence} 
shows how the algorithm behaves when we start with already non-collinear residuals.
This allows for the treatment of the case when the perturbed initial projection of the 
residual (\ref{eqn.shifted-incorrect-update}) renders collinearity invalid at the start.
\begin{theorem}\label{thm.gcrodr-shift-convergence}
	Suppose we begin the cycle as in {\rm (\ref{eqn.gcrodr-colin-resid})}, with approximate collinearity between the base and shifted residuals, satisfying the relation
\begin{equation}\label{eqn.initial-resids-noncollinear}
	\vek r_{0}^{(\sigma)} = \tilde{\beta}_{0}\vek r_{0} + \vek w^{(\sigma)}.
\end{equation}
If we perform a cycle of \RGMRES\  to reduce the residual of the base system and 
apply the approximate collinearity condition {\rm(\ref{eqn.approx-collinear-system})} to the shifted residual, 
then we have the relation
\begin{equation}\label{eqn.shifted-resid-relation2}
\tilde{\vek r}_{m}^{(\sigma)} = \tilde{\beta}_{m}\vek r_{m}  - \sigma \vek E\vek F\prn{\tilde{\vek y}_{m}^{(\sigma)}}_{1:k} + \vek w^{(\sigma)}.
\end{equation}
\end{theorem}
\textit{Proof.} We can write the residual produced by the approximate collinearity procedure for the shifted system
as follows, using (\ref{eqn.z-def}),

\begin{eqnarray*}
\tilde{\vek r}_{m}^{(\sigma)} & = & \vek b - (\vek A+\sigma\vek I)\vek x_{m}^{(\sigma)}\label{eqn.approx-collin-anal}\\
  & = &\vek r_{0}^{(\sigma)} - (\vek A+\sigma\vek I)\widehat{\vek V}_{m}\widetilde{\vek y}_{m}^{(\sigma)}\nonumber\\
  & = & \tilde{\beta}_{0}\vek r_{0} + \vek w^{(\sigma)} - (\vek A+\sigma\vek I)\widehat{\vek V}_{m}\tilde{\vek y}_{m}^{(\sigma)}\nonumber\\
  & = & \tilde{\beta}_{0}\vek r_{0} - \prn{\widehat{\vek W}_{m+1}\widetilde{\vek G}_{m}^{(\sigma)} + \sigma \brac{\begin{matrix}\vek E\vek F & \vek 0  \end{matrix}}}
\tilde{\vek y}_{m}^{(\sigma)} + \vek w^{(\sigma)} \nonumber\\
 & = & \tilde{\beta}_{0}\norm{\vek r_{0}}\widehat{\vek W}_{m+1}\vek e_{k+1}^{(m+1)}  - \widehat{\vek W}_{m+1}\widetilde{\vek G}_{m}^{(\sigma)}\tilde{\vek y}_{m}^{(\sigma)} - \sigma 
\brac{\begin{matrix}\vek E\vek F & \vek 0  \end{matrix}}\tilde{\vek y}_{m}^{(\sigma)} + \vek w^{(\sigma)}\nonumber\\
 & = & \tilde{\beta}_{0}\norm{\vek r_{0}}\widehat{\vek W}_{m+1}\vek e_{k+1}^{(m+1)}  - \widehat{\vek W}_{m+1}\widetilde{\vek G}_{m}^{(\sigma)}\tilde{\vek y}_{m}^{(\sigma)} - \tilde{\beta}
_{m}\widehat{\vek W}_{m+1}\vek z_{m+1} + \tilde{\beta}_{m}\widehat{\vek W}_{m+1}\vek z_{m+1} \nonumber\\ 
 &  & - \sigma \brac{\begin{matrix}\vek E\vek F & \vek 0  \end{matrix}}\tilde{\vek y}_{m}^{(\sigma)} + \vek w^{(\sigma)} \nonumber \\
 & = & \widehat{\vek W}_{m+1}\prn{\tilde{\beta}_{0}\norm{\vek r_{0}}\vek e_{k+1}^{(m+1)} - \widetilde{\vek G}_{m}^{(\sigma)}\tilde{\vek y}_{m}^{(\sigma)} -  \tilde{\beta}_{m}\vek z_{m+1}}+ 
\tilde{\beta}_{m}\widehat{\vek W}_{m+1}\vek z_{m+1}  \nonumber\\
 & & - \sigma \brac{\begin{matrix}\vek E\vek F & \vek 0  \end{matrix}}\tilde{\vek y}_{m}^{(\sigma)} + \vek w^{(\sigma)}.\nonumber
\end{eqnarray*}
Now using the approximate collinearity condition (\ref{eqn.approx-collinear-system}) and the fact that by definition $\vek r_{m} = \widehat{\vek W}_{m+1}\vek z_{m
+1}$, we have that 
\begin{equation} \nonumber \label{eqn.shifted-resid-rel}
	\tilde{\vek r}_{m}^{(\sigma)} = \tilde{\beta}_{m}\vek r_{m}  - \sigma \brac{\begin{matrix}\vek E\vek F & \vek 0  \end{matrix}}\tilde{\vek y}_{m}^{(\sigma)} + \vek w^{(\sigma)},
\end{equation}
which can be rewritten in the form (\ref{eqn.shifted-resid-relation2}).  \cvd

It should be noted that the term $-\sigma\vek E\vek F\prn{\tilde{\vek y}_{m}^{(\sigma)}}_{1:k} + \vek w^{(\sigma)}$ is a 
function of the quality of the recycled subspaces as well as of $\sigma$.  With the use of simple inequalities, we obtain 
an important corollary estimating the amount of residual norm reduction we can expect for the shifted systems.
\begin{corollary}\label{cor.shifted-resid-bound}
	The shifted system residual norm satisfies the following inequality,
\begin{equation}\label{eqn.shifted-resid-norm-estimate}
\norm{\tilde{\vek r}_{m}^{(\sigma)}} \leq \ab{\tilde{\beta}_{m}}\norm{\vek r_{m}}  + \ab{\sigma} \norm{\vek E\vek F}\norm{\prn{\tilde{\vek y}_{m}^{(\sigma)}}
_{1:k}} + \norm{\vek w^{(\sigma)}}.
\end{equation}
\end{corollary}

As long as 
$\ab{\tilde{\beta}_{m}}\norm{\vek r_{m}}$ dominates the right-hand side, we will observe 
a reduction of the shifted residual norm.
This reduction is controlled by $\ab{\sigma}
$, $\norm{\vek E\vek F}$, and $\norm{\prn{\tilde{\vek y}_{m}^{(\sigma)}}_{1:k}}$.  We cannot control $\norm{\prn{\tilde{\vek y}_{m}^{(\sigma)}}_{1:k}}$, and $\sigma$ is 
dictated by the problem.  The size of $\norm{\vek E\vek F}$ is connected to the quality of $\vek U$ as an approximation to an invariant subspace of $\vek A
$.  This can seen by writing 
\begin{equation}\label{eqn.EF}
	\vek E\vek F = \vek U - \prn{\vek C\vek Y + \vek V_{m+1}\vek Z}
\end{equation}
and observing that the norm of this difference decreases as $\vek U$ becomes a better approximation of an invariant subspace of $\vek A$.  
Thus, choosing $\vek U$ as an approximate invariant subspace may improve performance
of the method.

Ideally, we would like to detect
when $\ab{\tilde{\beta}_{m}}\norm{\vek r_{m}}$ ceases to dominate (\ref{eqn.shifted-resid-norm-estimate}) 
in order to cease updating the approximations to the shifted system
once such an update no longer leads to a decrease in residual norm. 
Our analysis gives us a way to monitor both quantities.   Observe that given $\sigma$, $\vek U$, and 
$\vek C$, if we compute $\widetilde{\vek y}_{m}^{(\sigma)}$ according to (\ref{eqn.approx-collinear-system}), 
then from (\ref{eqn.EF}), we can compute the product $\vek E\vek F\prn{\tilde{\vek y}_{m}^{(\sigma)}}_{1:k}$.  
Thus, we can keep track of the vector $\vek w^{(\sigma)}$, and use it to construct $\vek r_{m}^{(\sigma)}$ 
using (\ref{eqn.shifted-resid-relation2}).  Rather than detecting that 
$\ab{\tilde{\beta}_{m}}\norm{\vek r_{m}}$ ceases to dominate (\ref{eqn.shifted-resid-norm-estimate}),
it is simpler to calculate $\norm{\vek r_{m}^{(\sigma)}}$ after each cycle and detect when it has ceased
to be reduced by the correction from that cycle.  At this point, we cease updateing $\vek x^{(\sigma)}_{m}$
for the remaining cycles.

It should be noted that $\vek w^{(\sigma)}$ can be easily accumulated. 
At the beginning of Algorithm \ref{alg.GCRODR-approx-collinear}, we compute an initial value of $\vek w^{(\sigma)}$ 
according to (\ref{eqn.initial-resids-noncollinear}).  At Line \ref{alg.line.noproj} of Algorithm~\ref{alg.GCRODR-approx-collinear}, we 
update $\vek w^{(\sigma)} \leftarrow \vek w^{(\sigma)} - \sigma\vek U\vek C^{\ast}\widetilde{\vek r}^{(\sigma)}$ 
according to (\ref{eqn.proj-pert}).  At Line \ref{alg.line.approx-collin}, we 
update $\vek w^{(\sigma)}~\leftarrow~\vek w^{(\sigma)}~-~\sigma\vek E\vek F(\widetilde{\vek y}_{m})_{1:k}$ according to (\ref{eqn.shifted-resid-relation2}).

Does the linear system (\ref{eqn.approx-collinear-system}) correspond to an exact
collinear condition for some choice of deflation space?
Observe that if we write 
\begin{equation}\nonumber
	\vek U^{(\sigma)} = \vek U - \sigma(\vek A+\sigma\vek I)^{-1}\vek E\vek F,
\end{equation}
 then we obtain an exact Arnoldi-like relation  
\begin{equation}\label{eqn.new-arnoldi-rel}
	(\vek A+\sigma\vek I)\widehat{\vek V}_{m}^{(\sigma)} = \widehat{\vek W}_{m+1}\widetilde{\vek G}_{m}^{(\sigma)},
\end{equation}
where $\widehat{\vek V}_{m}^{(\sigma)} = \begin{bmatrix}  \vek U^{(\sigma)} & \vek V_{m-k}\end{bmatrix}$.
If we select $\vek x^{(\sigma)}_{m}\in\vek x_{0}^{(\sigma)}+\CR(\widehat{\vek V}_{m}^{(\sigma)})$ and enforce the collinearity condition $\vek r_{m}^{(\sigma)} = \beta_{m}\vek r_{m}$, 
then we see that (\ref{eqn.approx-collinear-system}) is the exact collinearity equation which must be solved to obtain the collinear residual.  Thus, the failure of the approximate collinearity condition,
due to singularity of (\ref{eqn.approx-collinear-system}), corresponds to the nonexistence
of an exactly collinear residual for the shifted system
over a different augmented subspace (which is unavailable in practice).


%% file: numerical_results.tex
\section{Numerical Experiments}\label{numresults}

\begin{figure}[htb]
\hfill
\includegraphics[scale=0.60]{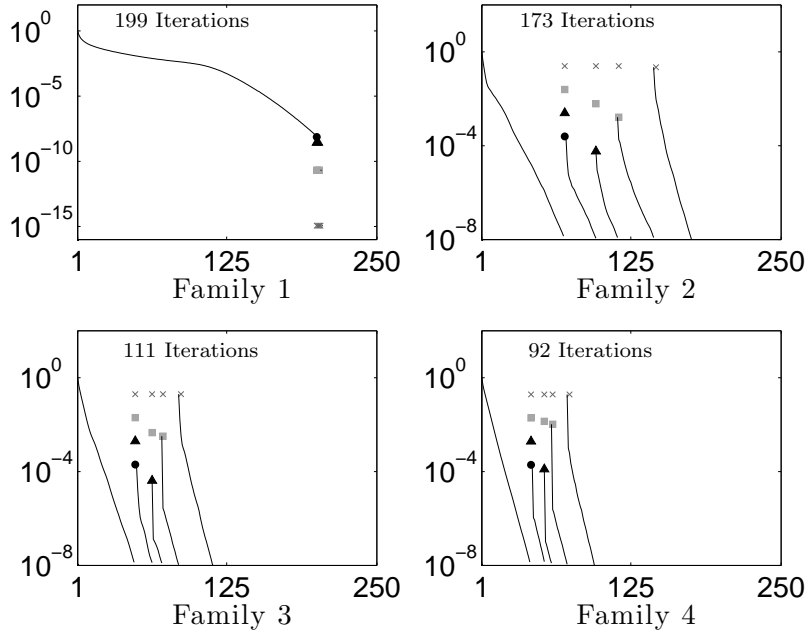}
\hfill
\begin{picture}(0,0)
\put(-140,230){{\footnotesize 173 Iterations}}
\end{picture}
\begin{picture}(0,0)
\put(-140,105){{\footnotesize 92 Iterations}}
\end{picture}
\begin{picture}(0,0)
\put(-300,230){{\footnotesize 199 Iterations}}
\end{picture}
\begin{picture}(0,0)
\put(-300,105){{\footnotesize 111 Iterations}}
\end{picture}
\caption{The performance of \RGMRES\ (RGMRES) for shifted systems with recursion on the number of unconverged shifted systems.  
We performed simple tests on a sequence 
\label{figure.shifted-recycled-bidiag-recursive}
of four families of $1000\times 1000$ bidiagonal matrices, each family with five systems.  
For this test, $m=100$ and $k=50$.  The first 
matrix is the bidiagonal matrix used in {\rm \cite{Darnell2008}}.  The 
other systems are constructed by applying $\CO(1)$ bidiagonal random perturbations to the first system using the {\ttfamily sprand()} Matlab function. The four shifts are 
$10^{-2}$, $10^{-1}$, $1$, and $10$.  For the shifted systems, the residuals are only computed at the end of each cycle,
 with residual norms represented by the \textbf{circle}, 
\textbf{triangle}, \textbf{square}, and \textbf{cross}, respectively.  The curves originating from these symbols are the convergence curves
for the \RGMRES\ iterations executed for each shift system when that system becomes the 
base system during a recursive call to the method.}
\end{figure}

We performed a series of numerical experiments illustrating the applicability of
the method described in Section~\ref{approx-collinear}, i.e., using an implementation
of Algorithm~\ref{alg.GCRODR-approx-collinear}.  
Following \cite{Parks.deSturler.GCRODR.2005},
we constructed recycled subspaces from harmonic Ritz vectors of the coefficient matrix associated with the base system~(\ref{eqn.Axb}) with respect 
to the augmented subspace. 
In the figures reported here, 
for each recursive call to the algorithm, the solid black line represents the convergence curve
for the base system, while the different markers indicate residual norms for the shifted systems
at the end of each restart cycle.

\begin{figure}[htb]
\hfill
\includegraphics[scale=0.60]{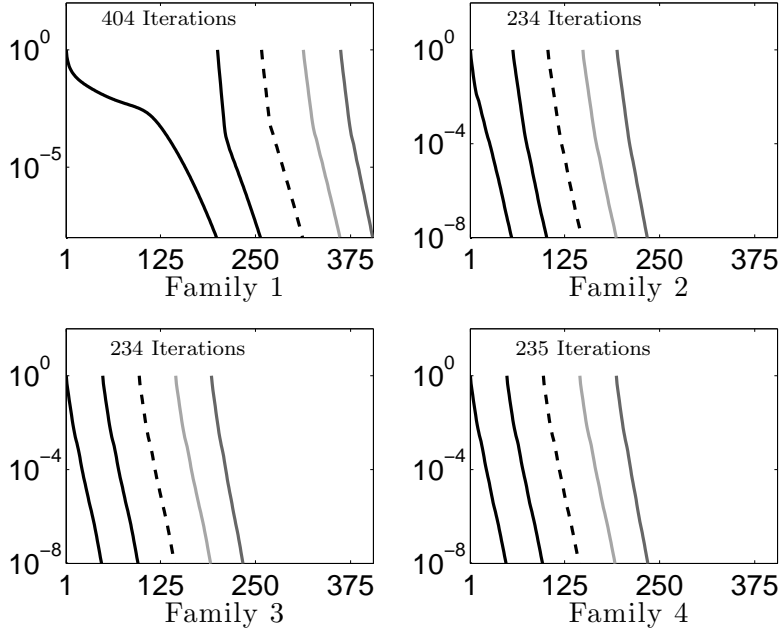}
\hfill
\begin{picture}(0,0)
\put(-140,230){{\footnotesize 234 Iterations}}
\end{picture}
\begin{picture}(0,0)
\put(-140,105){{\footnotesize 235 Iterations}}
\end{picture}
\begin{picture}(0,0)
\put(-300,230){{\footnotesize 404 Iterations}}
\end{picture}
\begin{picture}(0,0)
\put(-300,105){{\footnotesize 234 Iterations}}
\end{picture}
\caption{Convergence curves when we apply \RGMRES\ to each shifted bidiagonal system sequentially. As in the previous experiment, $m=100$ and $k=50$.
\label{figure.shifted-recycled-bidiag-gcrodr}}
\end{figure}  
\begin{figure}[htb]
\hfill
\includegraphics[scale=0.35]{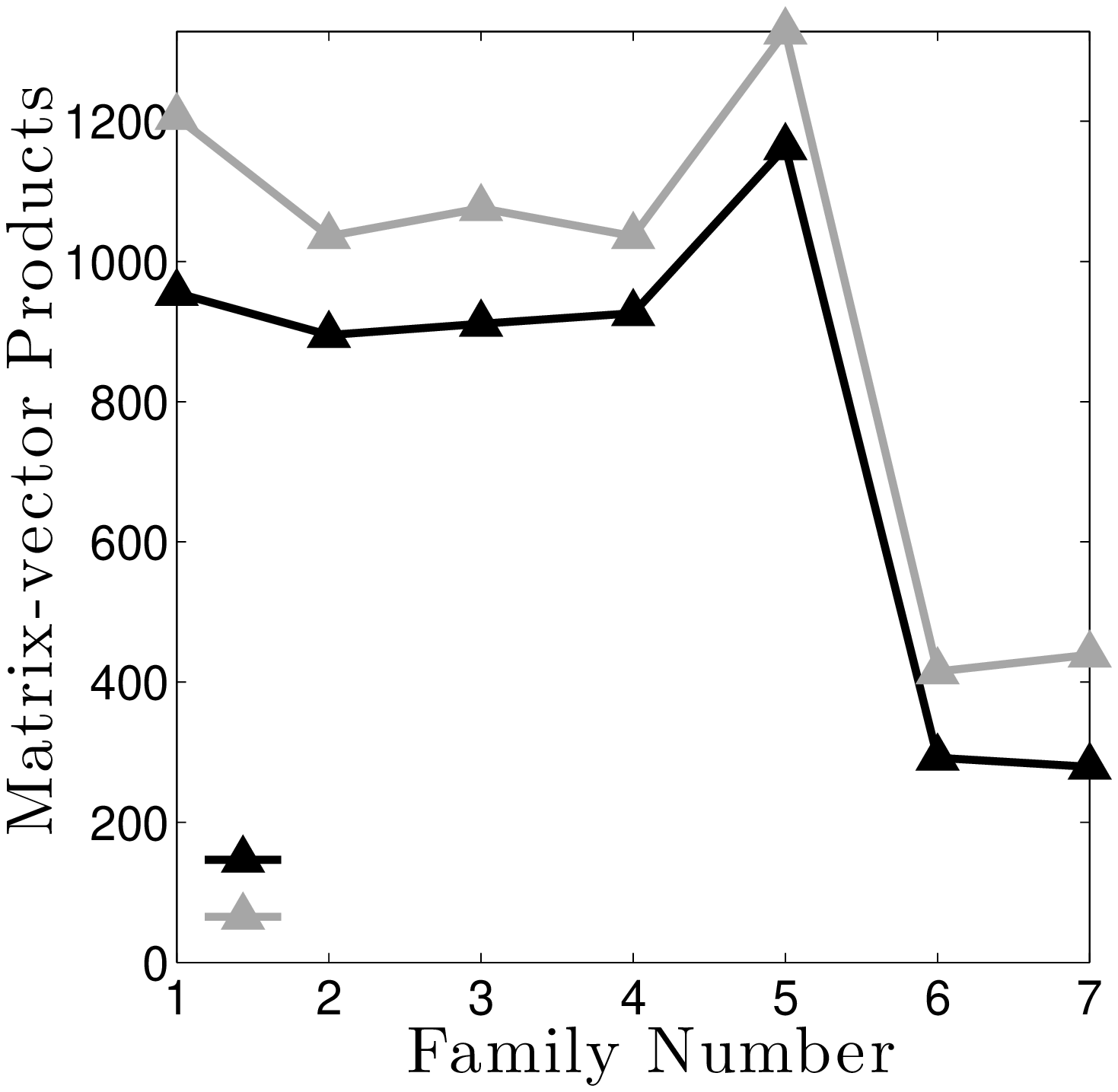}\quad\quad
\includegraphics[scale=0.35]{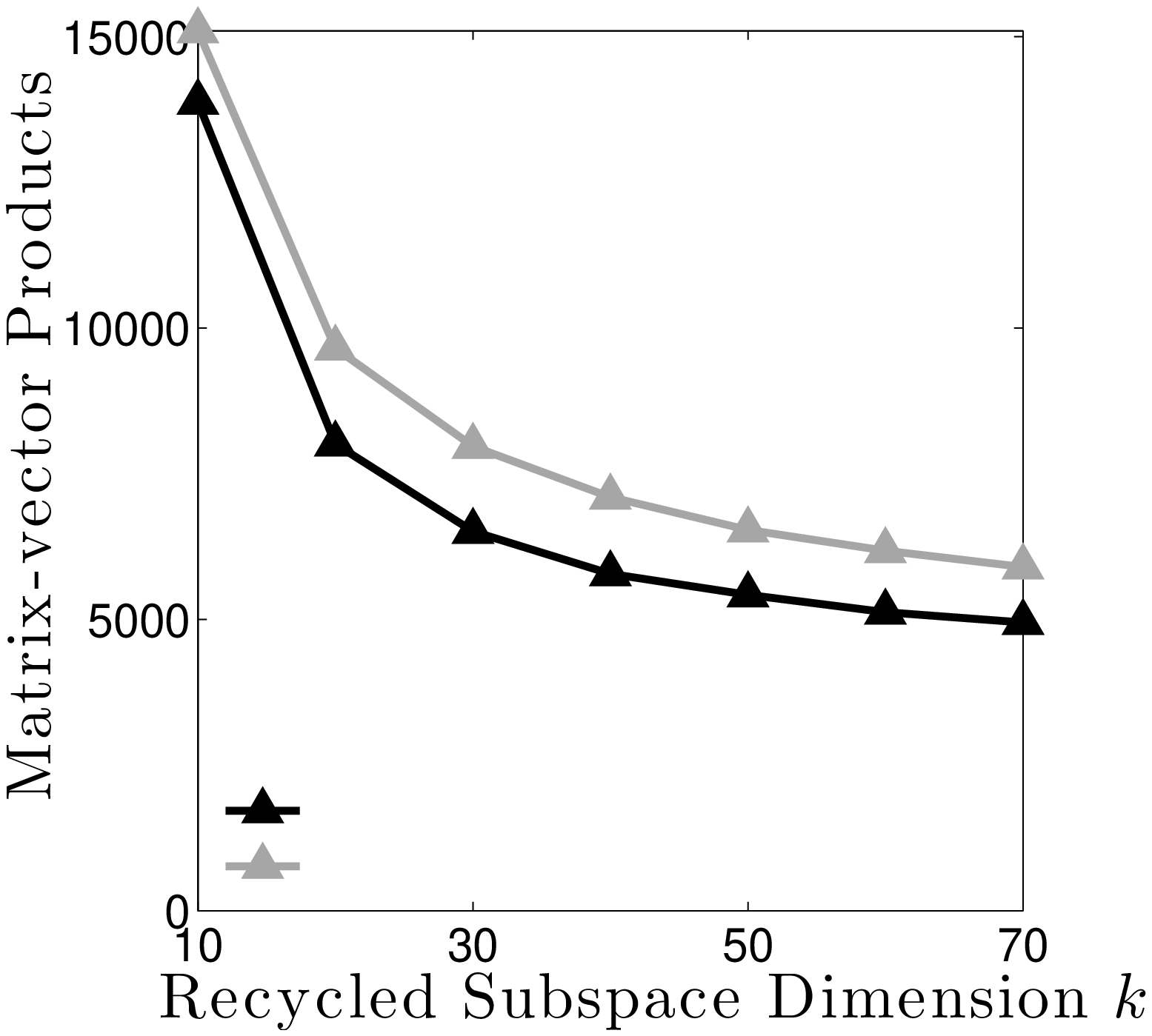}
\hfill
\begin{picture}(0,0)
\put(-117,30){\rotatebox{0}{\footnotesize Shifted RGMRES}}
\put(-117,20){\rotatebox{0}{\footnotesize Repeated RGMRES}}
\put(-328,30){\rotatebox{0}{\footnotesize Shifted RGMRES}}
\put(-328,20){\rotatebox{0}{\footnotesize Repeated RGMRES}}
\end{picture}
\caption{The performance of \RGMRES\ for shifted systems on a sequence of seven small Wilson fermion matrices where for each matrix, we solve a linear system with 
the base system and those associated to the shifts, $.001$ , $.002$, $.003$, $-.6$, and $-.5$.  In the left figure, we illustrate the performance of RMGRES($100$,$50$) for 
shifted systems as compared to repeated applications of RGMRES($100$,$50$).  In the figure on the right, we compare performance for different size recycled 
subspaces with $m=100$ fixed.}
\label{figure.shifted-recycle-qcd-small}
\end{figure}

In all experiments, 
when solving the first family of shifted systems in the sequence, there is no initial recycled subspace.  
Thus, (\ref{eqn.AU-UD.span}) holds at the end
of each cycle, and $\vek E=\vek 0$.  Therefore, the approximate collinearity condition 
(\ref{eqn.approx-collinear-system}) becomes an exact collinearity condition. 
This is equivalent to 
applying shifted GMRES-DR for shifted systems~\cite{Darnell2008}.  
Observe in the convergence plots, that all residuals 
are reduced in norm below tolerance when solving the base system.   

Our first experiment, presented in Figure \ref{figure.shifted-recycled-bidiag-recursive}, 
illustrates the performance of the \RGMRES\ method for shifted linear systems on a sequence of four 
bidiagonal matrices.  The first matrix, $\vek B_{1}$, used in \cite{Darnell2008}, is a bidiagonal matrix with  
$\curl{.1, 1, 2, \ldots, 998,999}$ on the diagonal 
and ones on the first superdiagonal, and the other matrices are random bidiagonal 
perturbations of $\vek B_{1}$, with the perturbations
having the same bidiagonal structure and having Frobenius norm $1$.  We see that, as predicted by Corollary 
\ref{thm.gcrodr-shift-convergence}, the amount of residual reduction achieved for the shifted systems is 
affected by the size of the shift.  For the shift $\sigma_{1} = 
10^{-2}$, the relative residual is reduced to $\CO(10^{-4})$ during the solution of the base system 
while for $\sigma_{4} = 10$, the relative residual is only reduced to $
\CO(10^{-1})$.  This experiment is more for illustrative purposes than to demonstrate superior performance.  
Nevertheless, after convergence for the base system, we take one of the shifted systems 
as our new base system and reapply the algorithm for the smaller family of systems.
\begin{figure}[htb!]
\hfill
\includegraphics[scale=0.60]{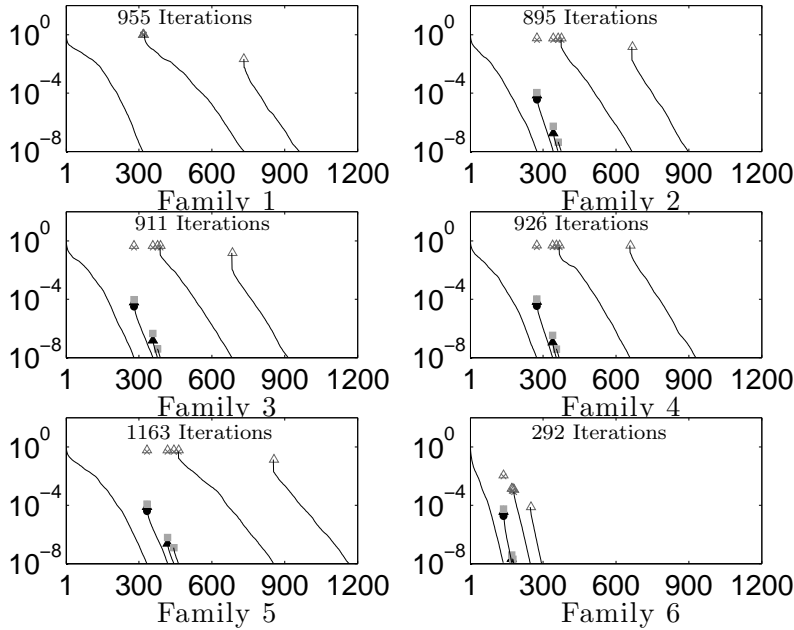}
\hfill
\begin{picture}(0,0)
\put(-280,152){{\footnotesize 911 Iterations}}
\end{picture}
\begin{picture}(0,0)
\put(-140,152){{\footnotesize 926 Iterations}}
\end{picture}
\begin{picture}(0,0)
\put(-140,230){{\footnotesize 895 Iterations}}
\end{picture}
\begin{picture}(0,0)
\put(-140,73){{\footnotesize 292 Iterations}}
\end{picture}
\begin{picture}(0,0)
\put(-300,230){{\footnotesize 955 Iterations}}
\end{picture}
\begin{picture}(0,0)
\put(-300,73){{\footnotesize 1163 Iterations}}
\end{picture}
\caption{Convergence curves from the same experiment as in Figure \ref{figure.shifted-recycle-qcd-small} 
but only for the first six systems. 
We again use $m=100$ and $k=50$.
 Observe that two of the shifted systems  (corresponding to the negative shifts) require more work 
 than the others for each base 
 system \textbf{including} the first system, in which we started with no recycled subspace.  
 These are the two shifts for which shifted 
 GMRES (or shifted GMRES-DR) would not converge.  Notice that in this case, we are still 
 able to converge by applying \RGMRES\ at the end.}
\label{figure.shifted-recycled-qcd-curves}
\end{figure}  
For comparison, we present in Figure~\ref{figure.shifted-recycled-bidiag-gcrodr} the 
convergence curves if we simply apply \RGMRES\ to each 
shifted system sequentially.
It can be appreciated that while applying \RGMRES\ sequentially requires a total of 1107
matrix-vector products for all systems with all shifts, the proposed approach requires
only 575 matrix-vector products, an improvement of about 50\%.

For our second and third experiments, we test two sequences of six
QCD matrices from the University of Florida sparse matrix collection 
\begin{figure}[htb!]
\begin{center}
\hfill
\includegraphics[scale=0.60]{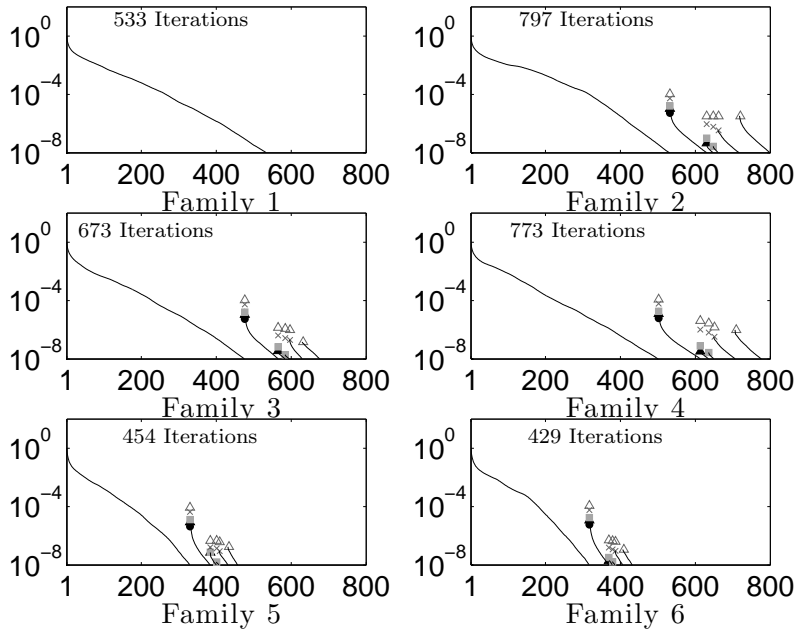}
\hfill
\begin{picture}(0,0)
\put(-300,150){{\footnotesize 673 Iterations}}
\end{picture}
\begin{picture}(0,0)
\put(-140,150){{\footnotesize 773 Iterations}}
\end{picture}
\begin{picture}(0,0)
\put(-140,230){{\footnotesize 797 Iterations}}
\end{picture}
\begin{picture}(0,0)
\put(-140,72){{\footnotesize 429 Iterations}}
\end{picture}
\begin{picture}(0,0)
\put(-300,230){{\footnotesize 533 Iterations}}
\end{picture}
\begin{picture}(0,0)
\put(-300,72){{\footnotesize 454 Iterations}}
\end{picture}
\vspace*{2em}
\end{center}
\caption{Convergence curves for another sequence of six Wilson fermion matrices, of 
size $49152\times 49152$ with $m=100$ and $k=50$.}
\label{figure.shifted-recycled-qcd-big-curves}
\end{figure} 
 \cite{DH.2011}.  In the second experiment, we work with six $3072\times 3072$ 
 sample matrices (called $\vek D_{1}$ through $\vek D_{6}$) with filename 
 prefix \verb|conf5.0-00l4x4|.  We can construct the coefficient matrix 
 $\vek A_{i} = \vek I -\kappa^{(i)} \vek D_{i}$ where $\kappa^{(i)}$ 
is a parameter associated to the QCD 
problem.  For each matrix, there exists some critical value $\kappa_{c}^{(i)}$ 
such that for  $0\leq \kappa^{(i)} < \kappa_{c}^{(i)}$, $\vek A_{i}$ is a real-positive 
matrix.  
Equivalently, for each $\vek A_{i}$, we can write $\vek A_{i} =  \frac{1}{\kappa^{(i)}}\vek I - \vek D_{i}$ 
where $\frac{1}{\kappa_{c}^{(i)}}< \frac{1}{\kappa^{(i)}}<\infty$, and we can scale any 
right-hand-side so that  we are solving the same problem. For each $\vek 
D_{i}$, ${\kappa_{c}^{(i)}}$ is included with the matrix, and in these experiments, all are 
in the interval $\brac{0.20,0.22}$.   Frequently in QCD computations, we wish to solve with multiple parameters.

We chose $\curl{.001, .002, .003, -.6, -.5}$ as our family of shifts.  Observe that by the 
definition of $\vek A_{i}$ and $\kappa_{c}^{(i)}$, the two shifted coefficient 
matrices associated with the two negative shifts are not real-positive.  These are not physically relevant for
QCD computations.  We chose negative shifts merely to demonstrate the robustness of the algorithm.
In this experiment, GMRES for shifted systems was unable to produce approximations 
for long sequences of iterations (due to numerical singularity of the augmented collinearity matrix).  
Since the shifted GMRES method did not converge for some systems, its performance was not included in the 
figure.  However, as we have noted, it is not difficult to modify this algorithm to gracefully 
handle this situation by applying restarted GMRES to any unconverged shifted 
systems at the end of the process.  We 
compared with another strategy, repeated applications of 
\RGMRES\ \cite{Parks.deSturler.GCRODR.2005} for the base and shifted 
system.  In Figure \ref{figure.shifted-recycle-qcd-small} we 
present the matrix-vector product counts for each system for a particular 
recycled subspace dimension as well as the totals over seven systems for various recycled 
subspace dimensions.  We see that our method is able to produce a $20\%$ reduction 
in the number of matrix-vector products needed to solve these systems, when compared 
to repeated applications of \RGMRES.  In Figure \ref{figure.shifted-recycled-qcd-curves}, 
we present the convergence curves for the first six QCD matrices.

In the third experiment, we worked with another sequence of six QCD matrices from 
the University of Florida Sparse Matrix Collection \cite{DH.2011} 
with filename prefixes \verb|conf5.4| and \verb|conf6.0|.  These matrices of size $49152\times 49152$.   
We used the critical $\kappa_{c}^{(i)}$ to construct our system matrices as in the second experiment, and
 we choose the shifts $\curl{.001, .002, .003, .01, .02}$ as in the second experiment.  
 In Figure \ref{figure.shifted-recycled-qcd-big-curves}, 
 we see the convergence of our algorithm for these systems.  

In the fourth experiment, we work with a sequence of eleven QCD matrices obtained 
from \cite{Wuppertal.Matrices.2012}.  These matrices were delivered
already shifted to be positive-real.  As in the previous experiment, they are also of size  $49152\times 49152$.   
\begin{table}[t!]
 \caption{ \label{table.qcd-defl-dim}
A comparison, in terms of iteration counts of the shifted GMRES algorithm (SGMRES)
with the \RGMRES\ algorithm for shifted systems (RGMRES) in terms of iteration count 
for different cycle lengths $m$.
The results presented are the total iterations for solving a sequence of eleven 
QCD systems from \cite{Wuppertal.Matrices.2012}.} 
\begin{center}
\begin{tabular}{r|c|c c c c } 
\multicolumn{1}{c|}{$m$}&$k$ & SGMRES($m$) &  RGMRES($m-k$,$k$)  & ratio\\
\hline
  25 &12& 4297 & 3880 & 0.90\\ 
   50&25& 3284 & 2980 & 0.91\\ 
  75 & 37 & 3108 & 2816& 0.91\\ 
  100 & 50 & 3028 & 2697 & 0.89\\ 
  125 & 67 & 3058 & 2612 &0.85\\ 
  150 & 75 &2958& 2546 &0.86\\ 
  175 & 87 & 2962 & 2499 &0.84\\ 
  200 & 100 &2947 & 2458 &0.83\\ 
  225 & 112 & 2860 & 2410 &0.84\\ 
\end{tabular} 
 \end{center}
 \end{table}
In Table \ref{table.qcd-defl-dim}, we illustrate the performance
of the proposed algorithm when the total dimension
of the augmented space increases. We also ran
a comparable instance of the shifted GMRES algorithm.
More specifically, 
we compared the performance of shifted GMRES with cycle length $m$ versus 
our algorithm with an $k=\floor{m/2}$ dimension deflation space and
$m-k$ cycle length. In this experiment, there are two shifts, $\curl{.8,.81}$.   
Here we see the potential benefits that our algorithm can yield as the deflation 
dimension increases.  
In particular, in the last example,
the gain as compared with RGMRES is of 16\%.
We mention though that for
a larger number of
shifts (or different values for the shifts) 
RGMRES may not be more advantagous.
This follows from the fact that
each shift incurs an additional recursive call to the algorithm and
additional iterations. There is no such increase for shifted GMRES. 
Therefore, for sufficiently large number 
of shifts, shifted GMRES will have an advantage.
Which method will perform better depends on several
factors including the number of shifts (as we just mentioned),
the magnitude of the shifts, the size of deflation space, and 
the deflation space selection technique. What we have shown is
that for certain problems, the recycling strategy is definitely
worth considering.

%% file: conclusions.tex
\section{Conclusions}\label{section.conclusions}
We have shown that 
the ideal method that solves a family of shifted systems simultaneously using one 
augmented subspace with subspace recycling generally does not exist under a fixed 
storage requirement independent of the number of shifts.
As an alternative, we present two methods, 
each of which relax one of the two requirements, yielding two possible algorithms.  One
solves the family of shifted systems using the same subspace but requires the construction of
multiple deflation spaces.  The other constructs approximations over a single augmented
Krylov subspace, but not all shifted systems are solved to tolerance at the same time.  
Instead, 
we showed some theoretical results indicating that the latter approach
produces improved initial approximate solutions for the shifted systems. This was
confirmed in our numerical experiments, which also showed
that the fixed-memory method can be quite effective, especially when
the shifts are all located in a small interval.



%% file: acknowledgement.tex
\section*{Acknowledgments}
 We would like to thank David Day and Michael Parks for engaging in fruitful discussions with the first author. We also thank 
 Andreas Frommer for his extensive comments, including the suggestion that Algorithm \ref{alg.GCRODR-approx-collinear} is amenable to recursion.